\documentclass[11pt]{article}

\usepackage{amssymb}
\usepackage{amsthm}
\usepackage{amsmath,amsfonts,amssymb}
\usepackage{graphicx}
\usepackage{algpseudocode}
\usepackage{algorithm}
\usepackage{natbib}
\usepackage{setspace}

\newtheorem{theorem}{Theorem}[section]

\newtheorem{corollary}[theorem]{Corollary}
\newtheorem{definition}[theorem]{Definition}
\newtheorem{remark}[theorem]{Remark}

\begin{document}

\title{Continuity of the Effective Path Delay Operator for Networks Based on the Link Delay Model}
\author{Ke Han$^{a}$
\qquad Terry L. Friesz$^{b}$\\\\
$^{a}$\textit{Department of Civil and Environmental Engineering}\\
\textit{Imperial College London, London SW7 2BU, UK}\\
\textit{Corresponding author. Email: k.han@imperial.ac.uk}\\\\
$^{b}$\textit{Department of  Industrial and Manufacturing Engineering}\\
\textit{Pennsylvania State University, University Park, PA 16802, USA}\\
\textit{Email: tfriesz@psu.edu}
}
\date{}
\maketitle

\begin{abstract}
This paper is concerned with a dynamic traffic network performance model, known as {\it dynamic network loading} (DNL), that is frequently employed in the modeling and computation of analytical {\it dynamic user equilibrium} (DUE). As a key component of continuous-time DUE models, DNL aims at describing and predicting the spatial-temporal evolution of traffic flows on a network that is consistent with established route and departure time choices of travelers, by introducing appropriate dynamics to flow propagation, flow conservation, and travel delays. The DNL procedure gives rise to the path delay operator, which associates a vector of path flows (path departure rates) with the corresponding path travel costs.  In this paper, we establish strong continuity of the path delay operator for networks whose arc flows are described by the {\it link delay model} \citep{Friesz1993}. Unlike result established in \cite{ZM}, our continuity proof is constructed 
 without assuming {\it a priori} uniform boundedness of the path flows. Such a more general continuity result has a few important implications to the existence of {\it simultaneous route-and-departure choice} DUE without {\it a priori} boundedness of path flows, and to any numerical algorithm that allows convergence to be rigorously analyzed. 
\end{abstract}

\section{Introduction}
{\it Dynamic traffic assignment} (DTA) is usually viewed as the descriptive modeling of time varying flows on vehicular networks consistent with established traffic flow. This paper is concerned with a specific type of dynamic traffic assignment known as continuous time simultaneous route-and-departure choice (SRDC) dynamic user equilibrium (DUE) for which unit travel cost, including early and late arrival penalties, is identical for those route and departure time choices selected by travelers between a given origin-destination pair. There are two essential components within the notion of DUE: (i) the mathematical expression of Nash-like equilibrium conditions, and (ii) a network performance model, which is, in effect, an embedded dynamic network loading (DNL) problem. The DNL relates closely to the {\it effective path delay operator}, which plays a pivotal role in any of the mathematical forms of DUE problem, including variational inequality \citep{Friesz1993}, differential variational inequality \citep{FBST, FM}, and complementarity problems \citep{Pang, Wie}.

Continuity of the effective delay operator, as we study in this paper, is critical to the DUE models since it is necessary to the existence of dynamic user equilibria \citep{Browder, existence, ZM}, and it is the minimum regularity assumption required by many computational algorithms that support convergence analysis, such as the fixed point algorithm \citep{FKKR, LWRDUE}, the descent direction method \citep{HL2003, SL2004}, and the projection method \citep{Ukkusuri}.

One way to model path delay in dynamic networks is that proposed by \cite%
{Friesz1993}, who employ arc exit time functions together with a
volume-dependent link traversal time function. Such a perspective on path
delay has been referred to both as the {\it link-delay model} (LDM)  and as the
{\it point queue model} (PQM) \citep{CTM}. Notice that despite the popular tendency among many scholars to refer to the Vickrey model \citep{Vickrey} as the point queue model, a careful literature search reveals that the name `point queue model' was first suggested by \cite{CTM} to describe the explicit travel time function model, in other words, the link delay model proposed by \cite{Friesz1993}. Since the \cite{Friesz1993} paper, conjectures but
few results about the qualitative properties of the LDM/PQM delay operator
have been published. One result that is needed for the study of  dynamic user
equilibrium existence as well as for analyses of DUE algorithms, when
network loading is based on the LDM/PQM, is continuity of the effective path delay
operator. It should be mentioned that \cite{ZM} investigated a similar problem and showed the weak continuity of the path delay operator under the assumption that the path flows are {\it a priori} bounded, where those bounds do not arise from any behavioral argument or theory. In particular, we note that {\it Theorem 5.1} of \cite{ZM}, which states the weak continuity of the delay operator, relies on the {\it strong first-in-first-out} (FIFO) assumption. Although the paper later stated that a sufficient condition for strong FIFO is uniform boundedness of path flows (path departure rates), it is not difficult to show that boundedness is also necessary to ensure strong FIFO. Effectively, the strong FIFO assumption is equivalent to the uniform boundedness of path flows.

Notably, this paper provides a more general continuity result, namely the path delay operator of interest is strongly continuous without the assumption of boundedness on path flows. Strong continuity without boundedness is central to the proof of SRDC DUE existence, and to the analyses of DUE algorithms. We point out that the simultaneous route-and-departure choice (SRDC) notion of DUE employs more general constraints relating path flows to a table of fixed trip volumes than the route choice (RC) DUE considered by \cite{ZM}. The boundedness assumption is less of an issue for the RC DUE by virtue of problem formulation: that is, for RC DUE, the travel demand constraints are of the following form:
\begin{equation}\label{demandconstraintrc}
\sum_{p\in\mathcal{P}_{ij}}h_p(t)~=~T_{ij}(t)\qquad\forall t,\quad\forall (i,\,j)\in\mathcal{W}
\end{equation}
where $\mathcal{W}$ is the set of origin-destination pairs, $\mathcal{P}_{ij}$ is the set of paths connecting $(i,\,j)\in\mathcal{W}$ and $h_p(t)$ is the departure rate along path $p$. Furthermore, $T_{ij}(t)$ represents the rate (not volume) at which travelers leave origin $i$ with the intent of reaching destination $j$ at time $t$; each such trip rate is assumed to be bounded from above. Since \eqref{demandconstraintrc} is imposed pointwise and every path flow $h_p$ is non-negative, we are assured that $h=(h_p:\,p\in\mathcal{P}_{ij},\,(i,\,j)\in\mathcal{W})$ are automatically uniformly bounded. On the other hand, the SRDC user equilibrium imposes the following constraints on path flows:
\begin{equation}\label{demandconstraintsrdc}
\sum_{p\in\mathcal{P}_{ij}}\int_{t_0}^{t_f}h_p(t)\,dt~=~Q_{ij}\qquad\forall (i,\,j)\in\mathcal{W}
\end{equation}
where $Q_{ij}\in\mathbb{R}_{++}$ is the volume (not rate) of travelers departing node $i$ with the intent of reaching node $j$. The integrals in \eqref{demandconstraintsrdc} are interpreted in the sense of Lebesgue; hence, \eqref{demandconstraintsrdc} alone is not enough to assure bounded path flows. This observation has been the major hurdle to providing existence without the invocation of bounds on path flows, and therefore, serves as the main motivation of our investigation of the continuity of the delay operator without assuming {\it a priori} boundedness on path flows.

\section{Background and preliminaries}

This paper is concerned with one type of dynamic traffic assignment (DTA) known as  simultaneous route-and-departure choice dynamic user equilibrium for which unit travel cost, including early and late arrival penalties, is identical for those route and departure time choices selected by travelers between a given originÐdestination pair. Such a model is originally presented in \cite{Friesz1993} and discussed subsequently by \cite{FBST, LWRDUE, FKKR}, and  \cite{FM}.

\subsection{Dynamic user equilibrium and the path effective delay operator}
We begin by  considering a planning horizon $[t_0,\, t_f]\subset \mathbb{R}_+$. Let $\mathcal{P}$ be the set of paths employed by road users. The most crucial ingredient of the DUE model  is the path delay operator. Such an operator, denoted by
\begin{equation*}
D_{p}(t,h)\text{ \ \ \ }\forall p\in \mathcal{P},
 \end{equation*}
maps  a given vector of departure rates $h$ to the collection of travel times. Each travel time is associated with a particular choice of route $p\in\mathcal{P}$ and departure time $t\in[t_0,\,t_f]$. The path delay operators usually do not take on any closed form, instead they can only be evaluated numerically through the dynamic network loading procedure. On top of the path delay operator we introduce the effective path delay operator which generalizes the notion of travel cost to include early or late arrival penalties. In this paper we consider the effective path delay operators of the following form.
\begin{equation}\label{psidef}
\Psi _{p}(t,h)=D_{p}(t,h)+F\left[ t+D_{p}(t,h)-T_{A}\right] \text{ \ \ \ }%
\forall p\in P
\end{equation}%
where $T_{A}$ is the target arrival time. In our formulation the target time $T_A$ is allowed to depend on the user classes. We introduce the fixed trip matrix $\big(Q_{ij}: (i,\,j)\in\mathcal{W}\big)$, where each $Q_{ij}\in \mathbb{R} _{+}$ is the fixed travel demand between origin-destination pair $\left( i,j\right) \in \mathcal{%
W}$. Note that $Q_{ij}$ represents traffic volume, not flow. Finally we let $\mathcal{P}_{ij}\subset \mathcal{P}$ to be the set of
paths connecting origin-destination pair $\left( i,j\right) \in \mathcal{W}
$. 

As mentioned earlier $h$ is the vector of path flows $h=\{h_p: p\in\mathcal{P}\}$. We stipulate that each path flow is square integrable, that is
$$
h\in\big(L_+^2[t_0,\,t_f]\big)^{|\mathcal{P}|}
$$
The set of feasible path flows is defined as
\begin{equation}\label{feasible}
\Lambda=\left\{ h\geq 0:\sum_{p\in \mathcal{P}_{ij}}%
\int_{t_{0}}^{t_{f}}h_{p}\left( t\right) dt=Q_{ij}\text{ \ \ \ }\forall
\left( i,j\right) \in \mathcal{W}\right\} \subseteq \left( L_{+}^{2}\left[
t_{0},t_{f}\right] \right) ^{\left\vert \mathcal{P}\right\vert }
\end{equation}
Let us also define the essential infimum of effective travel delays%
\begin{equation*}
v_{ij}=\hbox{essinf} \left[ \Psi _{p}(t,h):p\in \mathcal{P}_{ij}\right] \text{ \ \
\ \ }\forall \left( i,j\right) \in \mathcal{W}
\end{equation*}%
The following definition of dynamic user equilibrium was first articulated by \cite{Friesz1993}:
\begin{definition}
\label{duedef}{\bf (Dynamic user equilibrium)}. A vector of departure rates (path
flows) $h^{\ast }\in \Lambda $ is a dynamic user equilibrium if%
\begin{equation*}
h_{p}^{\ast }\left( t\right) >0,p\in P_{ij}\Longrightarrow \Psi _{p}\left[
t,h^{\ast }\left( t\right) \right] =v_{ij}
\end{equation*}%
We denote this equilibrium by $DUE\left( \Psi ,\Lambda,\left[
t_{0},t_{f}\right] \right) $.
\end{definition}
Using measure theoretic arguments, \cite{Friesz1993} established that
a dynamic user equilibrium is equivalent to the following variational
inequality under suitable regularity conditions:
\begin{equation}
\left. 
\begin{array}{c}
\text{find }h^{\ast }\in \Lambda _{0}\text{ such that} \\ 
\sum\limits_{p\in \mathcal{P}}\displaystyle \int\nolimits_{t_{0}}^{t_{f}}\Psi _{p}(t,h^{\ast
})(h_{p}-h_{p}^{\ast })dt\geq 0 \\ 
\forall h\in \Lambda %
\end{array}%
\right\} VI(\Psi ,\Lambda ,\left[ t_{0},t_{f}\right] )  \label{duevi}
\end{equation}%

\subsection{The link delay model}

We shall consider a general network $(\mathcal{A},\,\mathcal{V})$ where $%
\mathcal{A}$ and $\,\mathcal{V}$ denote the set of arcs and the set of
nodes, respectively. Additionally, we shall take the link traversal time to
be a linear function of the arc volume at the time of entry. We describe arc
volume as the sum of volumes associated with individual paths using this
arc, that is 
\begin{equation}
x_{a}(t)~=~\sum_{p\in \mathcal{P}}\delta _{ap}\,x_{a}^{p}(t)\qquad \forall
~a\in \mathcal{A}  \label{ldmdef1}
\end{equation}%
where $x_{a}^{p}(t)$ denotes the volume on arc $a\in \mathcal{A}$ associated
with path $p\in \mathcal{P}$, and $\mathcal{P}$ is set of all paths
considered. Each path $p$ $\in \mathcal{P}$ is taken to be the set of
consecutive arcs that consitute the path; that is 
\begin{equation*}
p=\{a_{1},\ldots ,a_{m(p)}\}
\end{equation*}
Furthermore we shall make use of the arc-path incidence matrix 
\begin{equation*}
\Delta =\left( \delta _{ap}\right) 
\end{equation*}%
where 
\begin{equation*}
\delta _{ap}~=~%
\begin{cases}
1\qquad  & \hbox{if arc}~a~\hbox{belongs to path}~p \\ 
0\qquad  & \hbox{otherwise}%
\end{cases}%
\end{equation*}%
We also let the time to traverse arc $a_{i}$ for drivers who arrive at its
entrance at time $t$ be denoted by $D_{a_{i}}\left( x_{a_{i}}\right) $. 

It will be convenient to introduce cumulative exit flows $V_{a_{i}}(\cdot
),\,V_{a_{i}}^{p}(\cdot )$:  
\begin{equation*}
V_{a_{i}}(t)~\doteq ~\int_{t_{0}}^{t}v_{a_{i}}(s)\,ds,\quad
V_{a_{i}}^{p}(t)~\doteq ~\int_{t_{0}}^{t}v_{a_{i}}^{p}(s)\,ds,\qquad p\in 
\mathcal{P},\,a_{i}\in p
\end{equation*}%
where the notation employed has obvious and conformal definitions relative
to that introduced previously. The following \textit{differential algebraic
equation} (DAE) system is another version of the DAE system employed by \cite%
{FKKR}: 
\begin{align}
\displaystyle X_{a_{i}}^{p}(t)& ~=~V_{a_{i-1}}^{p}(t)-V_{a_{i}}^{p}(t)
\label{ldmdef5} \\
\displaystyle V_{a_{i-1}}^{p}(t)& ~=~V_{a_{i}}^{p}\big(%
t+D_{a_{i}}(X_{a_{i}}(t))\big),\qquad i~=~1,\ldots ,m(p)
\end{align}%
Furthermore, we make the following definitions 
\begin{equation}
v_{a_{0}}^{p}(t)~\equiv ~h_{p}(t),\qquad V_{a_{0}}^{p}(t)~\equiv
~\int_{t_{0}}^{t}h_{p}(s)\,ds  \label{ldmdef6}
\end{equation}%
where $h_{p}$, also known as the path flow, is the rate of departure from
the origin of $p\in \mathcal{P}$. It is also conventional to introduce the
link exit time function 
\begin{equation*}
\tau _{a}(t)\equiv t+D_{a}\left[ x_{a}(t)\right] 
\end{equation*}
for each $a\in \mathcal{A}$. 

The next theorem, from \cite{Friesz1993}, presents an important property of
linear link delay functions:

\begin{theorem}
\label{ldmthm} For any linear arc delay function of the form $D(x)=\alpha
x+\beta,\, \alpha,\,\beta>0$, the resulting arc exit time function $\tau$ is
continuous and strictly increasing and hence the inverse function $\tau^{-1}$
exists.
\end{theorem}

\begin{proof}
See Theorem 1 of \cite{Friesz1993}.
\end{proof}

\section{The main result}

The following is a statement of our main result:

\begin{theorem}
\label{LDMlemma} Consider a general network $(\mathcal{A},\,\mathcal{V})$,
where arc dynamics are governed by the link delay model, assume the link
delay function for each $a\in \mathcal{A}$ is affine. That is 
\begin{equation*}
D_{a}\left[ x_{a}(t)\right] ~=~\alpha _{a}X_{a}(t)+\beta _{a}
\end{equation*}
where $\alpha _{a}\in \Re _{+}^{1}$ and $\beta _{a}\in \Re _{++}^{1}$. Then
the effective delay operator from $\Lambda $ into $\big(L^{2}[t_{0},\,t_{f}]\big)^{|\mathcal{P}|}:h\in \Lambda \longrightarrow \Psi
(\cdot ,\,h)$ is a continuous map.
\end{theorem}
\begin{remark}
In \cite{ZM}, the authors show that the effective delay operator is weakly continuous when the LDM is employed, under the restrictive assumption that the path flows are a priori bounded from above. Such assumption is dropped in our result; we also assert strong, not weak continuity.
\end{remark}

\noindent The continuity result for the effective delay operator from Theorem \ref{LDMlemma} is crucial for the analysis and computation of dynamic user equilibrium, especially when {\it a priori} upper bound on the path flows is not guaranteed by any behavioral or mathematical arguments.

\subsection{Proof of the main result}

Now we present the proof of Theorem \ref{LDMlemma}. 

\begin{proof}

We begin by showing that given a converging sequence $h^{(n)}$ in the space $%
\Lambda \subset \Big(L_{+}^{2}([t_{0},\,t_{f}])\Big)^{|\mathcal{P}%
|}$ such that 
\begin{equation}
\big\|h^{(n)}-h\big\|_{L^{2}}~\longrightarrow ~0\qquad
n~\longrightarrow ~\infty ,  \label{ldml2conv}
\end{equation}%
the corresponding delay function $D_{p}\big(\cdot ,\,h^{(n)}\big)$ converges
uniformly to $D_{p}\big(\cdot ,\,h\big)$ for all $p\in \mathcal{P}$. This
will be proved in several steps.

\noindent \textbf{Part} \textbf{1.} First, let us consider just one single
arc, and hence omit the subscript $a$ for brevity. Assume a sequence of
entering flows $\{u^{(n)}\}_{n\geq 1}$ converging to $u$ in the $L%
^{2}([t_{0},\,t_{f}])$ space; that is 
\begin{equation}
\big\|u^{(n)}-u\big\|_{2}~\doteq ~\left( \int_{t_{0}}^{t_{f}}\Big(%
u^{(n)}(t)-u(t)\Big)^{2}dt\right) ^{1/2}~\longrightarrow ~0\qquad
n~\longrightarrow ~\infty   \label{l2conv}
\end{equation}%
Define the cumulative entering vehicle counts 
\begin{equation*}
\begin{array}{l}
\displaystyle U^{(n)}(t)~\doteq ~\int_{t_{0}}^{t}u^{(n)}(s)\,ds\qquad n~\geq
~1 \\ 
\displaystyle U(t)~\doteq ~\int_{t_{0}}^{t}u(s)\,ds%
\end{array}%
\qquad t\in \lbrack t_{0},\,t_{f}]
\end{equation*}%
Then we assert the uniform convergence $U^{(n)}\longrightarrow U$ on $%
[t_{0},\,t_{f}]$. To see this, fix any $\varepsilon >0$, in view of (\ref%
{l2conv}), choose $N>0$ such that for all $n>N$ 
\begin{equation*}
\big\|u^{(n)}-u\big\|_{2}~<~\varepsilon 
\end{equation*}%
According to the embedding of $L^{1}([t_{0},\,t_{f}])$ into $%
L^2[t_0,\,t_f]$, we deduce for any $t\in \lbrack
t_{0},\,t_{f}]$ that 
\begin{align*}
\big|U^{(n)}(t)-U(t)\big|& ~=~\left\vert
\int_{t_{0}}^{t}u^{(n)}(s)\,ds-\int_{t_{0}}^{t}u(s)\,ds\right\vert  \\
& \\
& ~\leq ~\big\|u^{(n)}-u\big\|_{1}~\leq ~(t_{0}-t_{f})^{1/2}\,\big\|u^{(n)}-u%
\big\|_{2} \\
& \\
& ~<~(t_{0}-t_{f})^{1/2}\varepsilon 
\end{align*}%
The preceding shows the uniform convergence $U^{(n)}(\cdot)\longrightarrow U(\cdot)$ on $%
[t_{0},\,t_{f}]$.

\noindent \textbf{Part 2.} We adapt the recursive technique devised in \cite%
{Friesz1993}. Let $X^{(n)}(\cdot )$, \thinspace $n\geq 1$, and $X(\cdot )$
denote the arc volumes corresponding to $U^{(n)}(\cdot ),\,n\geq 1$ and $%
U(\cdot )$, respectively. Assume, without loss of generality, that 
\begin{equation*}
X^{(n)}(t_{0})~=~0,\qquad X(t_{0})~=~0
\end{equation*}%
and that, for the flow profile $U(\cdot )$, the first vehicle enters the arc
of interest at time $t_{0}$. In addition, let $t_{1}$ denote the time that
first vehicle exits the arc of interest. By definition 
\begin{equation}
t_{1}~=~D(0)~=~\beta   \label{ldmeqn1}
\end{equation}%
For all $t\in \lbrack t_{0},\,t_{1}]$, since no vehicle can exit the arc
before time $t_{1}$, we have 
\begin{equation}
X^{(n)}(t)~=~U^{(n)}(t)\qquad X(t)~=~U(t)\qquad t\in \lbrack t_{0},\,t_{1}]
\label{ldmeqn2}
\end{equation}%
For each flow profile $U^{(n)},\,n\geq 1$, denote the exit time function
restricted to $[t_{0},\,t_{1}]$ by $\tau _{1}^{(n)}(\cdot )$. Under the flow
profile $U$, denote the exit time function restricted to $[t_{0},\,t_{1}]$
by $\tau _{1}(\cdot )$. Then 
\begin{align}
\tau _{1}^{(n)}(t)& ~=~t+D\big(X^{(n)}(t)\big)~=~t+a\,U^{(n)}(t)+\beta
,\qquad t\in \lbrack t_{0},\,t_{1}]  \label{ldmeqn3} \\
\tau _{1}(t)& ~=~t+D\big(X(t)\big)~=~t+a\,U(t)+\beta ,\qquad t\in \lbrack
t_{0},\,t_{1}]
\end{align}%
We conclude that $\tau _{1}^{(n)}\longrightarrow \tau _{1}$ uniformly on $%
[t_{0},\,t_{1}]$. Now let 
\begin{equation*}
\tilde{t}_{2}~\doteq ~\inf_{n}\tau _{1}^{(n)}(t_{1})~\leq ~\tau _{1}(t_{1})
\end{equation*}

\begin{figure}[htbp]
\centering
\includegraphics[width=0.6\textwidth]{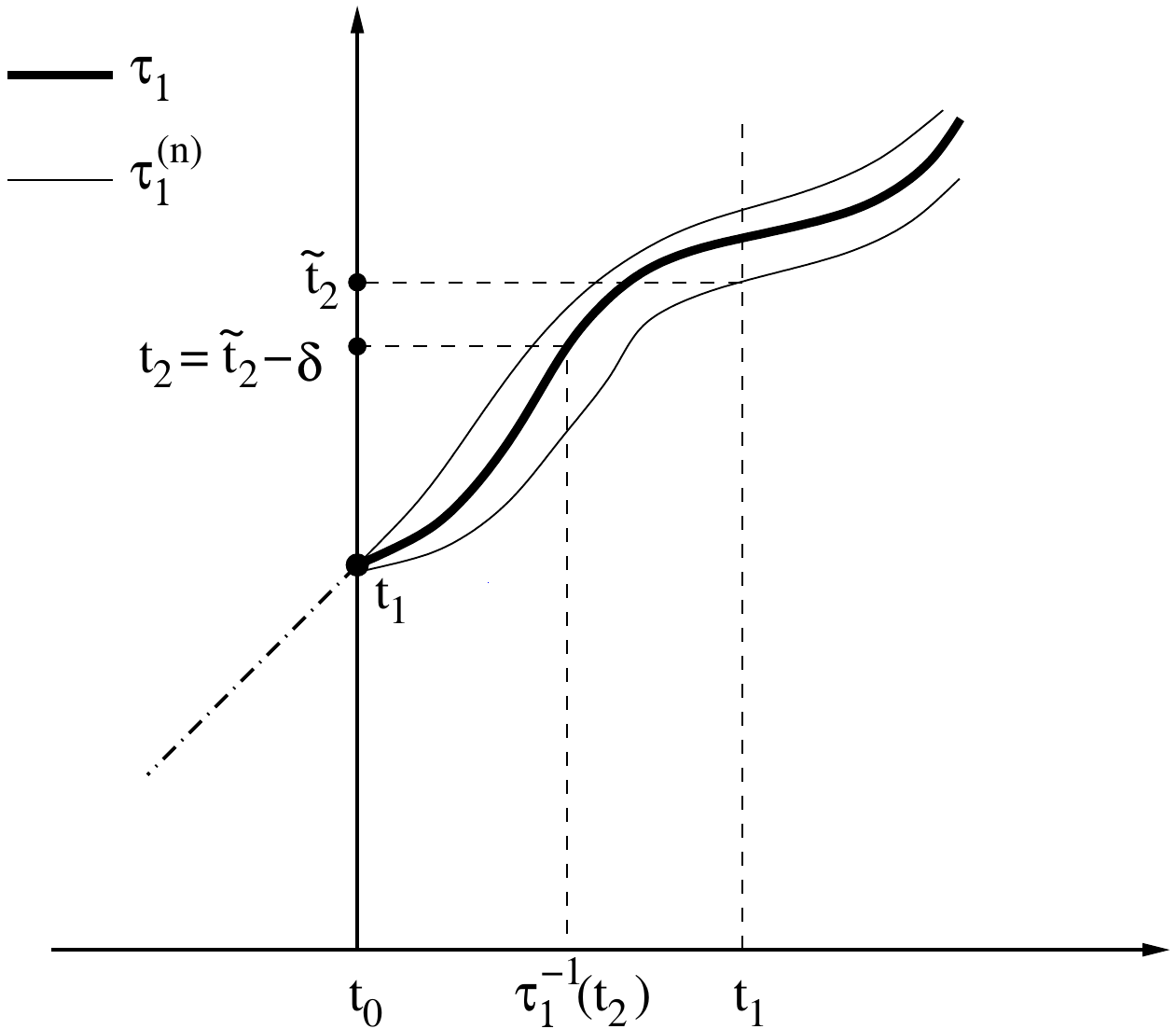}
\caption{Definitions of $\tilde t_2$ and $t_2$.}
\label{figproof1}
\end{figure}

\noindent Fix $\delta $ small enough and call $t_{2}\doteq \tilde{t}_{2}-\delta $.  See Figure \ref{figproof1} for a graphical illustration of these notations. By Theorem \ref{ldmthm} the functions $\big(\tau _{1}^{(n)}\big)^{-1}$, \thinspace $n\geq 1$%
, and $\tau _{1}^{-1}$, are well-defined, continuous and strictly increasing.
We claim that $\{\big(\tau _{1}^{(n)}\big)^{-1}\}_{n\geq 1}$ uniformly
converges to $\tau _{1}^{-1}$ on $[t_{1},\,t_{2}]$. To see this, we need to
extend the arrival time function $\tau _{1}$ and $\tau _{1}^{(n)}$ to the
interval $(-\infty ,\,t_{0})$. Because no vehicle is present during $%
(-\infty ,\,t_{0})$, it is natural to assign 
\begin{equation*}
\tau _{1}(t)~=~t+\beta ,\qquad \tau ^{(n)}(t)~=~t+\beta 
\end{equation*}%
This means, if an infinitesimal flow particle enters the arc at $t\in
(-\infty ,\,t_{0})$, its travel delay will always be $\beta $. Fix any $%
\varepsilon <t_{1}-\tau _{1}^{-1}(t_{2})$, and consider the following
quantities: 
\begin{equation}
\Delta _{\varepsilon }^{-}~\doteq ~\inf_{t\in \lbrack t_{0},\,\tau
_{1}^{-1}(t_{2})]}\Big\{\tau _{1}(t)-\tau _{1}(t-\varepsilon )\Big\}\qquad
\Delta _{\varepsilon }^{+}~\doteq ~\inf_{t\in \lbrack t_{0},\,\tau
_{1}^{-1}(t_{2})]}\Big\{\tau _{1}(t+\varepsilon )-\tau _{1}(t)\Big\}
\label{nb}
\end{equation}%
Since the infimum of a continuous function on a compact interval must be
obtained at some point $t\in \lbrack t_{0},\,\tau _{1}^{-1}(t_{2})]$, we
conclude $\Delta _{\varepsilon }^{-},\,\Delta _{\varepsilon }^{+}>0$ by the
strict monotonicity of $\tau _{1}$ established in Theorem \ref{ldmthm}.

According to the uniform convergence $\tau _{1}^{(n)}\longrightarrow \tau
_{1}$ on $(-\infty ,\,t_{1}]$, there exists some $N>0$ such that as soon as $%
n\geq N$, we have 
\begin{equation}
\tau _{1}^{(n)}(t)~\leq ~\tau _{1}(t)+\Delta _{\varepsilon }^{-}/2\qquad \
\tau _{1}^{(n)}(t)~\geq ~\tau _{1}(t)-\Delta _{\varepsilon }^{+}/2\qquad \
\forall t\in \lbrack -\infty ,\,t_{1}]  \label{unifconv}
\end{equation}%
For any $s\in \lbrack t_{1},\,t_{2}]$, we have $\tau _{1}^{-1}(s)\in \lbrack
t_{0},\,\tau _{1}^{-1}(t_{2})]$. Therefor, for all $n\geq N$, in view of (%
\ref{unifconv}) and (\ref{nb}), we have 
\begin{align*}
& \tau _{1}^{(n)}\big(\tau _{1}^{-1}(s)-\varepsilon \big)~\leq ~\tau _{1}%
\big(\tau _{1}^{-1}(s)-\varepsilon \big)+\Delta _{\varepsilon }^{-}/2~\leq
~\tau _{1}\big(\tau _{1}^{-1}(s)\big)-\Delta _{\varepsilon
}^{-}/2~=~s-\Delta _{\varepsilon }^{-}/2 \\
& \tau _{1}^{(n)}\big(\tau _{1}^{-1}(s)+\varepsilon \big)~\geq ~\tau _{1}%
\big(\tau _{1}^{-1}(s)+\varepsilon \big)-\Delta _{\varepsilon }^{+}/2~\geq
~\tau _{1}\big(\tau _{1}^{-1}(s)\big)+\Delta _{\varepsilon
}^{+}/2~=~s+\Delta _{\varepsilon }^{+}/2
\end{align*}%
By the Intermediate Value Theorem, there exists some $t^{\ast }\in \lbrack
\tau _{1}^{-1}(s)-\varepsilon ,\,\tau ^{-1}(s)+\varepsilon ]$ with $\tau
_{1}^{(n)}(t^{\ast })=s$. In other words, we know 
\begin{equation*}
\big|\big(\tau _{1}^{(n)}\big)^{-1}(s)-\tau _{1}^{-1}(s)\big|~=~\big|t^{\ast
}-\tau _{1}^{-1}(s)\big|~<~\varepsilon ,\qquad \forall ~n~\geq ~N
\end{equation*}%
This finishes our claim that $\big(\tau _{1}^{(n)}\big)^{-1}\longrightarrow
\tau _{1}^{-1}$ uniformly on $[t_{1},\,t_{2}]$.

Let $\tau _{2}^{(n)}(\cdot ),\,n\geq 1$ and $\tau _{2}(\cdot )$ be the exit
time functions for commuters entering during the interval $[t_{1},\,t_{2}]$,
corresponding to entering flow profiles $U^{(n)}(\cdot ),\,n\geq 1$ and $%
U(\cdot )$, respectively. Then for each $t\in \lbrack t_{1},\,t_{2}]$, we
may state 
\begin{align}
\tau _{2}^{(n)}(t)& ~=~t+a\Big[U^{(n)}(t)-U^{(n)}\big((\tau
_{1}^{(n)})^{-1}(t)\big)\Big]+\beta ,\quad n~\geq ~1  \label{ldmeqn4} \\
\tau _{2}(t)& ~=~t+a\Big[U(t)-U\big(\tau _{1}^{-1}(t)\big)\Big]+\beta  \label{ldmeqn5}
\end{align}%
Now we make the claim that $U^{(n)}\Big(\big(\tau _{1}^{(n)}\big)^{-1}(t)%
\Big)\longrightarrow U\big(\tau _{1}^{-1}(t)\big)$ uniformly on $%
[t_{1},\,t_{2}]$. Indeed, for any $\varepsilon >0$, there exists $N_{1}>0$
such that, if $n>N_{1}$, we have 
\begin{equation*}
\big|U^{(n)}(t)-U(t)\big|~<~\varepsilon /2,\qquad \forall ~t\in \lbrack
t_{0},\,t_{f}]
\end{equation*}%
Moreover, the functions $\big(\tau _{1}^{(n)}\big)^{-1}$,\thinspace $\ n\geq
1$, and $\tau _{1}^{-1}$ have a uniformly bounded range on $[t_{1},\,t_{2}]$%
, namely $[t_{0},\,t_{1}]$. By the Heine-Cantor theorem, $U(\cdot )$
restricted to $[t_{0},\,t_{1}]$ is uniformly continuous, which means we can
find $\delta _{0}>0$ such that, for any $s_{1},\,s_{2}\in \lbrack
t_{0},\,t_{1}]$ with $|s_{1}-s_{2}|<\delta _{0}$, the following holds: 
\begin{equation*}
\big|U(s_{1})-U(s_{2})\big|~<~\varepsilon /2
\end{equation*}%
By uniform convergence of $\big(\tau _{1}^{(n)}\big)^{-1}\longrightarrow
\tau _{1}^{-1}$, we may choose $N_{2}>0$ so that, for $n>N_{2}$, we have 
\begin{equation*}
\big|\big(\tau _{1}^{(n)}\big)^{-1}(t)-\tau _{1}^{-1}(t)\big|~<~\delta _{0}
\end{equation*}%
Thus we deduce that, given $n>\max \{N_{1},\,N_{2}\}$, for any $t\in \lbrack
t_{1},\,t_{2}]$, the following is true: 
\begin{align*}
& \Big|U^{(n)}\Big(\big(\tau _{1}^{(n)}\big)^{-1}(t)\Big)-U\Big(\tau
_{1}^{-1}(t)\Big)\Big| \\
~\leq ~& \Big|U^{(n)}\Big(\big(\tau _{1}^{(n)}\big)^{-1}(t)\Big)-U\Big(\big(%
\tau _{1}^{(n)}\big)^{-1}(t)\Big)\Big|+\Big|U\Big(\big(\tau _{1}^{(n)}\big)%
^{-1}(t)\Big)-U\Big(\tau _{1}^{-1}(t)\Big)\Big| \\
~<~& \varepsilon /2+\varepsilon /2~=~\varepsilon 
\end{align*}%
This shows the uniform convergence $U^{(n)}\Big(\big(\tau _{1}^{(n)}\big)%
^{-1}(t)\Big)\longrightarrow U\big(\tau _{1}^{-1}(t)\big)$ on $%
[t_{1},\,t_{2}]$, and our claim is substantiated. As an immediate
consequence of (\ref{ldmeqn4}) and (\ref{ldmeqn5}), $\tau _{2}^{(n)}$
converges to $\tau _{2}$ uniformly on $[t_{1},\,t_{2}]$.

\noindent \textbf{Part 3.} We now proceed by induction as follows. Choose
any $\nu \geq 2$, and call 
\begin{equation*}
\tilde{t}_{\nu +1}~\doteq ~\inf_{n}\tau _{\nu }^{(n)}(t_{\nu }),\qquad
t_{\nu +1}~\doteq ~\tilde{t}_{\nu +1}-\delta 
\end{equation*}%
where the constant $\delta $ is the same as what was used in \textbf{Part 2}%
. Using the induction hypothesis that $\tau _{\nu }^{(n)}$ converges to $%
\tau _{\nu }$ uniformly on $[t_{\nu -1},\,t_{\nu }]$, we show the following
uniform convergence on $[t_{\nu },\,t_{\nu +1}]$: 
\begin{equation*}
\big(\tau _{\nu }^{(n)}\big)^{-1}~\longrightarrow ~\tau _{\nu }^{-1},
\end{equation*}%
The proof is similar to what has been done in \textbf{Part 2}. Now introduce 
$\tau _{\nu +1}^{(n)}(\cdot )$, $\,n\geq 1$, and $\tau _{\nu +1}(\cdot )$,
which are the exit time functions corresponding to $U^{(n)}(\cdot ),\,n\geq 1
$ and $U(\cdot )$ respectively; and they are restricted to the time interval 
$[t_{\nu },\,t_{\nu +1}]$. Similar to results (\ref{ldmeqn4}) and (\ref%
{ldmeqn5}), we have for $t\in \lbrack t_{\nu },\,t_{\nu +1}]$, that the
following holds: 
\begin{align}
\tau _{\nu +1}^{(n)}(t)& ~=~t+a\Big[U^{(n)}(t)-U^{(n)}\big(\big(\tau _{\nu
}^{(n)}\big)^{-1}(t)\big)\Big]+\beta \text{ }\quad n~\geq ~1  \label{ldmthm7}
\\
\tau _{\nu +1}(t)& ~=~t+a\Big[U(t)-U\big(\tau _{\nu }^{-1}(t)\big)\Big]%
+\beta 
\end{align}%
It can be shown as before that $U^{(n)}\big(\big(\tau _{\nu }^{(n)}\big)^{-1}%
\big)\longrightarrow U\big(\tau _{\nu }^{-1}\big)$ uniformly on $[t_{\nu
},\,t_{\nu +1}]$. Therefore $\tau _{\nu +1}^{(n)}\longrightarrow \tau _{\nu
+1}$ uniformly on $[t_{\nu },\,t_{\nu +1}]$. This finishes the induction.

\noindent \textbf{Part }\noindent \textbf{4.} We now have obtained a sequence $\big\{%
\lbrack t_{\nu },\,t_{\nu +1}]\big\}_{\nu \geq 0}$ of intervals. On each
interval $[t_{\nu },\,t_{\nu +1}]$, the uniform convergence 
\begin{equation*}
\tau _{\nu +1}^{(n)}~\longrightarrow ~\tau _{\nu +1},
\end{equation*}%
holds. Notice that, by construction, $t_{\nu +1}-t_{\nu }\geq \beta -\delta
>0$, $\forall \nu \geq 1$; therefore the interval $[t_{0},\,t_{f}]$ can be
covered by finitely many such intervals. As a consequence, we easily obtain
the uniform convergence of $\tau ^{(n)}(\cdot )\longrightarrow \tau (\cdot )$
on the whole of $[t_{0},\,t_{f}]$, where $\tau ^{(n)}(\cdot )$ corresponds
to $U^{(n)}(\cdot ),\,n\geq 1$ and $\tau (\cdot )$ corresponds to $U(\cdot )$%
.

Let $v^{(n)}(\cdot )$, \thinspace $n\geq 1$, and $v(\cdot )$ be the exit
flows of the single arc and then define the cumulative exit vehicle count 
\begin{equation*}
V(t)~\doteq ~\int_{t_{0}}^{t}v(s)\,ds,\qquad V^{(n)}(t)~\doteq
~\int_{t_{0}}^{t}v^{(n)}(s)\,ds,
\end{equation*}%
It then follows immediately from the relationships 
\begin{equation*}
V(t)~=~U\big(\tau ^{-1}(t)\big)\qquad V^{(n)}(t)~=~U^{(n)}\Big(\big(\tau
^{(n)}\big)^{-1}(t)\Big)
\end{equation*}%
that $V^{(n)}$ converges uniformly to $V(\cdot )$ on $[t_{0},\,t_{f}]$.

\noindent \textbf{Part 5.} Consider a general network $\big(\mathcal{A},\,%
\mathcal{V}\big)$ with a converging sequence $h^{(n)}\longrightarrow h$ in $%
\Lambda \subset \big(L^2[t_0,\,t_f]\big)^{|\mathcal{P}|}$.
Define for $p\in \mathcal{P}$ the following: 
\begin{equation*}
H_{p}^{(n)}(t)~\equiv ~\int_{t_{0}}^{t}h_{p}^{(n)}(s),\,ds,\qquad
H_{p}(t)~\doteq ~\int_{t_{0}}^{t}h_{p}(s)\,ds
\end{equation*}%
Then the $H_{p}^{(n)}(\cdot )$ converge uniformly to $H_{p}(\cdot )$. For
each arc $a\in \mathcal{A}$, where the notation employed has an obvious
meaning, the cumulative entering vehicle count $U_{a}^{(n)}(\cdot )$ is
given by 
\begin{equation*}
U_{a}^{(n)}(t)~=~\sum_{p}H_{p}^{(n)}(t)+\sum_{a^{\prime }\in \mathcal{I}%
(a)}V_{a^{\prime }}^{(n)}(t)
\end{equation*}%
In the above, the first summation is over all paths that use $a$ as the
first arc; and, in the second summation, $\mathcal{I}(a)$ denotes the set of
arcs immediately upstream from arc $a$. A simple mathematical induction with
the results established in previous steps yields the uniform convergence 
\begin{equation*}
U_{a}^{(n)}(t)~\longrightarrow ~U_{a}(t),\quad
V_{a}^{(n)}(t)~\longrightarrow ~V_{a}(t),\quad \tau
_{a}^{(n)}(t)~\longrightarrow ~\tau _{a}(t),\quad \forall a\in \mathcal{A}
\end{equation*}%
where $\tau _{a}(\cdot )$ is the exit time function of arc $a$. Thus, for
each path $p\in \mathcal{P}$, the path delay $D_{p}\big(\cdot ,\,h^{(n)}\big)
$ also converges uniformly to $D_{p}\big(\cdot ,\,h\big)$ since it is a
finite sum of arc delays. It remains to show that the effective delays obey $%
\Psi _{p}\big(\cdot ,\,h^{(n)}\big)\longrightarrow \Psi _{p}\big(\cdot
,\,h^{(n)}\big)$ uniformly on $[t_{0},\,t_{f}]$. Recall 
\begin{equation*}
\Psi _{p}(t,\,h)~=~D_{p}(t,\,h)+\mathcal{F}\big(t+D_{p}(t,\,h)-T_{A}\big)
\end{equation*}%
Notice that $\mathcal{F}$ is uniformly continuous on $[t_{0},\,t_{f}]$ by
the Heine-Cantor theorem; this means, for any $\varepsilon >0$, there exists 
$\sigma >0$ such that whenever $|s_{1}-s_{2}|<\sigma $, we have 
\begin{equation*}
\big|\mathcal{F}(s_{1})-\mathcal{F}(s_{2})\big|~<~\varepsilon /2
\end{equation*}%
Moreover, by uniform convergence, there exits $N>0$ such that, for all $n>N$%
, we have 
\begin{equation*}
\big|D_{p}\big(t,\,h^{(n)}\big)-D_{p}\big(t,\,h\big)\big|~<~\min \{\sigma
,\varepsilon /2\}\qquad \forall ~t\in \lbrack t_{0},\,t_{f}]
\end{equation*}%
We then readily deduce that, given $n>N$, the following holds for all $t\in
\lbrack t_{0},\,t_{f}]$:  
\begin{align*}
& \Big|\Psi _{p}\big(t,\,h^{(n)}\big)-\Psi _{p}\big(t,\,h\big)\Big| \\
~\leq ~& \Big|D_{p}\big(t,\,h^{(n)}\big)-D_{p}\big(t,\,h\big)\Big|+\Big|%
\mathcal{F}\Big(t+D_{p}\big(t,\,h^{(n)}\big)-T_{A}\Big)-\mathcal{F}\Big(%
t+D_{p}\big(t,\,h\big)-T_{A}\Big)\Big| \\
~<~& \varepsilon /2+\varepsilon /2~=~\varepsilon 
\end{align*}

\noindent \textbf{Part 6}. Finally, uniform convergence on the compact
interval $[t_{0},\,t_{f}]$ implies convergence in the $L^{2}$
norm: 
\begin{equation}
\int_{t_{0}}^{t_{f}}\Big(\Psi _{p}\big(t,\,h^{(n)}\big)-\Psi _{p}\big(t,\,h%
\big)\Big)^{2}\,dt~\longrightarrow ~0,\quad n~\longrightarrow ~\infty
,\qquad p\in \mathcal{P}  \label{ldmeqn8}
\end{equation}%
Summing up (\ref{ldmeqn8}) over $p\in \mathcal{P}$, we obtain the
convergence $\Psi \big(\cdot ,\,h^{(n)}\big)\longrightarrow \Psi \big(\cdot
,\,h\big)$ in the Hilbert space $\big(L^2[t_0,\,t_f]\big)
^{|\mathcal{P}|}$. 
\end{proof}

In some analysis \citep{BH, E-DUE} a slightly different notion of continuity of the effective delay operator is invoked as follows. 

\begin{definition}{\bf (A-continuity)} We say that the effective path delay operator $\Psi:\Lambda\to \big(L_2[t_0,\,t_f]\big)^{|\mathcal{P}|}$ is A-continuous if, for any weakly convergent sequence $h^{(n)}\rightharpoonup h^*\in\Lambda$ such that $|h_p^{(n)}(t)|<C$ for any $t\in[t_0,\,t_f]$ and $p\in\mathcal{P}$ where $C>0$ is some fixed constant, the effective path delays $\Psi_p(\cdot,\,h^{(n)})$ converges uniformly to $\Psi_p(\cdot,\,h^*)$ for each $p\in\mathcal{P}$. 
\end{definition}

As the following corollary asserts, the A-continuity also holds for the effective path delay operator. 

\begin{corollary}
Under the same assumption made in Theorem \ref{LDMlemma}, the effective delay operator $\Psi:\Lambda\to \big(L^2[t_0,\,t_f]\big)^{|\mathcal{P}|}$ is A-continuous.
\end{corollary}
\begin{proof}
Given any weakly convergent sequence $h^{(n)}\rightharpoonup h^*\in\Lambda$ such that $|h_p^{(n)}(t)|<C$ for any $t\in[t_0,\,t_f]$ and $p\in\mathcal{P}$, we define 
$$
H^{(n)}(t)~\doteq~\int_{t_0}^t h^{(n)}(s)\,ds,\qquad H^*(t)~\doteq~\int_{t_0}^t h^*(s)\,ds\qquad\forall t\in[t_0,\,t_f]
$$
Then one immediately has that $H^{(n)}(\cdot)$ converges to $H^*(\cdot)$ uniformly on $[t_0,\,t_f]$. Then the rest of the proof is the same as {\bf Part 2} - {\bf Part 6} of the proof of Theorem \ref{LDMlemma}.
\end{proof}

\section{Concluding remarks}

Dynamic traffic assignment differs from static traffic assignment in that
path delay does not enjoy a closed form. In fact the path delays needed for
the study of dynamic user equilibria (DUE) are operators that may only be
specified numerically. Moreover, such path delay operators are based on the
specific model of arc delay employed for network loading. We have considered
in this paper path delay operators for the network loading procedure that is
endogenous to Friesz et al. (1993) which has been referred to as the
link delay model (LDM) and also as the point queue model (PQM). We have
shown that LDM/PQM path delay operators are strongly continuous under the very mild
assumption that the link traversal time function is affine. In addition, our proof of continuity relies on no {\it ad hoc} assumptions on the uniform boundedness of path flows.  Combined with the results, \cite{Browder, HFY1, HFY2, existence}, an increasingly clear understanding of DUE based on different network performance models is emerging.

\end{document}